\newif\ifviewchanges
\viewchangesfalse
\documentclass{elsart}
\usepackage[sumlimits]{amsmath}
\usepackage{amsfonts}
\usepackage{ulem}
\usepackage{graphicx,graphics,color}
\DeclareMathOperator{\diff}{d}
\newcommand{\pp}[2]{\frac{\partial #1}{\partial #2}} 

\newcommand{\dd}[2]{\frac{\diff#1}{\diff#2}}

\newcommand{\Ro}{\mathrm{Ro}}
\newcommand{\Fr}{\mathrm{Fr}}
\def\MM#1{\boldsymbol{#1}}

\newcommand{\iso}{{\mathrm{iso}}}
\newcommand{\PP}{{\mathrm{P}}}
\newcommand{\q}{{\mathrm{q}}}

\usepackage{xspace}
\newcommand{\pdgp}{$\textrm{P1}_{\textrm{DG}}$-$\textrm{P2}$\xspace}
\newcommand{\pndgpn}{$\textrm{Pn}_{\textrm{DG}}$-$\textrm{P(n+1)}$\xspace}

\usepackage{natbib}
\ifviewchanges
\definecolor{fill}{rgb}{0.0,0.0,1.0}

\newcommand{\remove}[1]{\sout{#1}}
\newcommand{\mathremove}[1]{\mbox{\sout{$#1$}}}
\definecolor{dhfill}{rgb}{0.6,0.0,0.6}

\newcommand{\dhremove}[1]{\sout{#1}}
\newcommand{\dhmathremove}[1]{\mbox{\sout{$#1$}}}
\else
\newcommand{\remove}[1]{}
\newcommand{\mvec}[1]{\MM{#1}}
\newcommand{\mathremove}[1]{}

\newcommand{\eqnref}[1]{(\ref{#1})}

\newcommand{\dhremove}[1]{}
\newcommand{\dhmathremove}[1]{}
\fi

\begin{document}
\begin{frontmatter}
\title{A mixed discontinuous/continuous finite element pair for 
shallow-water ocean modelling}
\author[aeroIC]{Colin J. Cotter\corauthref{cor}}
\corauth[cor]{Corresponding author.}
\ead{colin.cotter@imperial.ac.uk}
\author[eseIC]{David A. Ham}
\ead{d.ham@imperial.ac.uk}
\author[eseIC]{Christopher C. Pain}
\ead{c.pain@imperial.ac.uk}
\address[aeroIC]{ Department of Aeronautics,
Imperial College London, London SW7 2AZ, United Kingdom}
\address[eseIC]{Department of Earth Science and Engineering,
Imperial College London, London SW7 2AZ, United Kingdom}
\begin{abstract}
  We introduce a mixed discontinuous/continuous finite element pair
  for ocean modelling, with continuous quadratic layer thickness and
  discontinuous velocity. We investigate the finite element pair
  applied to the linear shallow-water equations on an $f$-plane. The
  element pair has the property that all geostrophically balanced
  states which strongly satisfy the boundary conditions have discrete
  divergence equal to exactly zero and hence are exactly steady states
  of the discretised equations. This means that the finite element
  pair has excellent geostrophic balance properties. We also show that
  the element pair applied to the non-rotating linear shallow-water
  equations does not have any spurious small eigenvalues. We
  illustrate these properties using numerical tests and provide
  convergence calculations which show that the numerical solutions
  have errors which decay quadratically with element edge length for
  both velocity and layer thickness.
\end{abstract}

\end{frontmatter}

\section{Introduction}
A number of finite element pairs have been proposed for the rotating
shallow-water equations, including the $\PP1_{NC}-\PP1$ and
$\PP1-\iso$ $\PP2 - P1$ elements (investigated and compared to
several other element pairs in \citep{Ro_etal1998}), the RT0 elements
(introduced in \cite{RaTh1977} and proposed for the shallow-water
equations in \cite{WaCa1998}) and equal-order elements with
stablisation (also proposed in \cite{WaCa1998}); all of these elements
have been shown to perform well when integrating the rotating
shallow-water equations. In this paper we investigate the numerical
properties of the \pdgp finite element pair applied to the linear
shallow-water equations on an $f$-plane in order to investigate the
suitability of the element for shallow-water ocean modelling. The
finite element pair consists of discontinuous linear elements for
velocity and continuous quadratic elements for layer thickness. Even
though the layer thickness has shape functions which are one order higher
than velocity, there are still more degrees of freedom in the space of
discontinuous linear functions than the space of continuous quadratic
functions (except in meshes with very few elements), which is a
necessary (but not sufficient) condition for the absence of spurious
pressure modes (modes with high spatial frequency but small
eigenvalues in the wave operator which can pollute the numerical
solution with noise). In \citet{CoHaPaRe2007}, it was shown that in
one dimension, this choice leads to a discretisation of the wave
equation without rotation which does not have any spurious modes. It
was also shown that the dispersion relation is very accurate for the
first half of the discrete spectrum. A number of numerical tests were
carried out on two- and three-dimensional unstructured meshes which
showed that there were no spurious eigenvalues present. In fact, as we
show in this paper, the discretised Laplacian obtained by combining
the first-order discrete divergence and gradient operators is the same
as the usual Galerkin finite element discretised Laplacian obtained by
multiplying the Laplacian by a test function and integrating by
parts. Consequentially, the discretised equations without rotation do
not have any spurious eigenvalues, also the element pair applied to
the incompressible Navier-Stokes equations leads to an LBB-stable
discretisation without spurious pressure modes. For a general
discussion of LBB stability, see \citet{GrSa2000,KaSh2005}.
\citet{KaSh2005} also contains an exposition of the discontinuous
Galerkin method. For applications of the discontinuous Galerkin method
to waves equations see \citet{AiMoMu2006}, and for some applications
of the discontinuous Galerkin method to the rotating shallow-water
equations see \citet{AmBo2007,Levin2006,Be_etal2007,Gi2006}.

In this paper we concentrate on the interaction of the geostrophic
modes with the inertia-gravity waves which is crucial to the good
representation of large-scale dynamics. We find that not only does the
finite element pair allow for accurate representation of
geostrophically-balanced states, these states are completely uncoupled
from the inertia-gravity waves: the states are exactly steady as in
the unapproximated partial-differential equations. In section
\ref{mixed} we introduce the element pair applied to the linear
shallow water equations, show that the element pair has a discrete
Laplacian without spurious eigenvalues, and show that the element pair
has exactly steady geostrophic states. In section \ref{numerical
  tests} we verify these results with numerical tests. We also show
numerical calculations using Kelvin waves which are geostrophically
balanced in one direction; these waves are a good test of preservation
of balance. The results do not show any radiating inertia-gravity
waves. We provide convergence test results using the Kelvin wave exact
solution which confirm that the errors spatial discretisation
converges quadratically for both velocity and layer thickness,
indicating that the element pair is stable. Finally we give a summary
and outlook in section \ref{summary}.

\section{The mixed element}
\label{mixed}
In this section we describe our mixed element formulation applied to
the linear shallow-water equations on an $f$-plane.
\subsection{Mixed continuous/discontinuous Galerkin discretisation}
We start with the linearised 
shallow-water equation on an $f$-plane in
non-dimensional units
\begin{eqnarray}
\label{wave 1}
\mvec{u}_t + \frac{1}{\Ro}\mvec{k}\times\mvec{u} 
+ \frac{1}{\Fr^2}\nabla h &=& 0, \quad \mvec{u} = (u_1,\ldots,u_d), \\
\quad h_t + \nabla\cdot\mvec{u} &=& 0, \label{wave 2}
\end{eqnarray}
where $\mvec{u}$ is the velocity, $h$ is the perturbation layer
thickness, $\mvec{k}$ is the unit vector in the $z$-direction,
$\Ro=U/fL$ is the Rossby number, $\Fr=\sqrt{U/gH}$ is the Froude
number, $U$ is a velocity scale, $L$ is a horizontal length scale, $H$
is the mean layer thickness, $f$ is the Coriolis parameter and $g$ is
the acceleration due to gravity. The boundary conditions are
\begin{equation}
\label{bcs}
\mvec{u}\cdot\mvec{n} = 0 \quad\mathrm{on}\quad \partial\Omega
\end{equation}
where $\partial\Omega$ denotes the boundary of the domain $\Omega$, and
$\mvec{n}$ is the normal to $\partial\Omega$. To obtain the
discontinuous/continuous Galerkin form of the equations we multiply equation
(\ref{wave 1}) by a discontinuous test function $\mvec{w}$ and equation
(\ref{wave 2}) by a continuous test function $\phi$ and integrate over an
element $E$ to obtain
\begin{eqnarray}
\label{weak 1}
\dd{}{t}\int_{E} \mvec{w}\cdot\mvec{u}\diff{V} 
+ \frac{1}{\Ro}\int_{E} \mvec{w}\cdot\mvec{k}\times\mvec{u}\diff{V}
& = &
-\frac{1}{\Fr^2}\int_{E}\mvec{w}\cdot\nabla h\diff{V}, \\
\dd{}{t}\int_{E} \phi h \diff{V} & = & 
-\int_{E} \phi\nabla\cdot\mvec{u} \diff{V}.
\label{weak 2}
\end{eqnarray}
We then integrate equation (\ref{weak 2}) by parts, and make use of
the boundary conditions \eqnref{bcs} to obtain 
\begin{eqnarray}
\dd{}{t}\int_{E} \mvec{w}\cdot\mvec{u}\diff{V} 
+ \frac{1}{\Ro}\int_{E} \mvec{w}\cdot\mvec{k}\times\mvec{u}\diff{V}
& = &
-\frac{1}{\Fr^2}
\int_{E}\mvec{w}\cdot\nabla h\diff{V}, \\
\dd{}{t}\int_{E} \phi h \diff{V} & = & 
\int_{E} \nabla\phi\cdot\mvec{u} \diff{V} \\
& & \qquad - \int_{\partial E\backslash\partial\Omega}
\mvec{n}\cdot\tilde{\mvec{u}}\phi \diff{S},
\end{eqnarray}
where $\tilde{\mvec{u}}$ is the value of $\mvec{u}$ on the element
boundary $\partial E$, determined by the particular choice of
discontinuous Galerkin scheme which is chosen (the value on the upwind
face, for example), and where $\MM{n}$ is the outward-pointing unit
normal to the surface $\partial E$. Conservation requires that
$\tilde{\mvec{u}}$ takes the same value on either side of each
face. We sum these equations over all elements and the surface terms
cancel since $\phi$ is continuous.  This gives the form of the
equations that we will discretise:
\begin{eqnarray}
\label{weak parts 1}
\dd{}{t}\int_{\Omega} \mvec{w}\cdot\mvec{u}\diff{V} 
+ \frac{1}{\Ro}\int_{\Omega} \mvec{w}\cdot\mvec{k}\times\mvec{u}\diff{V}
& = &
-\frac{1}{\Fr^2}
\int_{\Omega}\mvec{w}\cdot\nabla h\diff{V}, \\
\label{weak parts 2}
\dd{}{t}\int_{\Omega} \phi h \diff{V} & = & 
\int_{\Omega} \nabla\phi\cdot\mvec{u} \diff{V}.
\end{eqnarray}
Derivatives are only applied to the scalar functions $h$ and $\phi$
and not the vector functions $\mvec{u}$ and $\mvec{w}$ which we shall
discretise with discontinuous elements. To add nonlinear advection it
is necessary to develop surface integrals on the boundaries of the
elements following the standard discontinuous Galerkin finite element
approach.

\subsection{The \pdgp element}
\begin{figure}
\begin{center}
\includegraphics*[width=4cm]{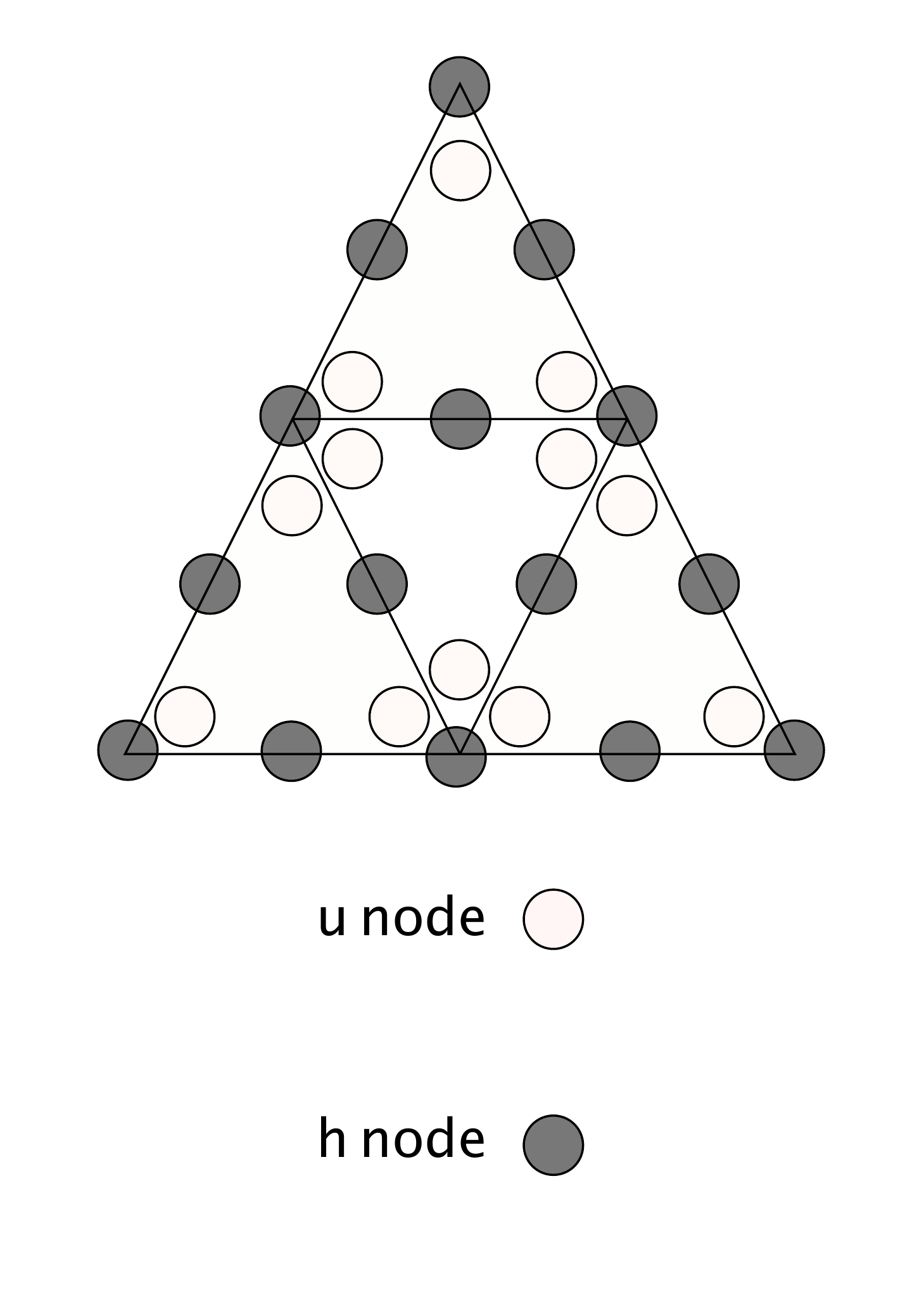}
\end{center}
\caption{\label{dof} Figure showing the distribution of nodes in the
  two-dimensional \pdgp element. Each element contains three nodes for each
  of the two components $u$ and $v$ of velocity, and six nodes for the
  layer thickness, but the latter nodes are shared across element boundaries
  since the layer thickness space is continuous.}
\end{figure}
In this subsection we develop the \pdgp discretisation for the
shallow-water equations. We make the choice that $\mvec{u}$ and
$\mvec{w}$ are approximated by discontinuous linear finite element
functions $\mvec{u}^\delta$ and $\mvec{w}^{\delta}$, whilst $\phi$ and
$h$ are approximated by continuous quadratic linear finite element
functions $h^\delta$ and $\phi^\delta$.

The Galerkin finite element approximation of equations (\ref{weak
  parts 1},\ref{weak parts 2}) is then
\begin{eqnarray*}
\dd{}{t}\int_{\Omega} \mvec{w}^\delta\cdot\mvec{u}^\delta\diff{V} 
+ \frac{1}{\Ro}\int_{\Omega} \mvec{w}^\delta\cdot\mvec{k}\times\mvec{u}^\delta
\diff{V}
& = &
-\frac{1}{\Fr^2}
\int_{\Omega}\mvec{w}^\delta\cdot\nabla h^\delta\diff{V}, \\
\dd{}{t}\int_{\Omega} \phi^\delta h^\delta \diff{V} & = & 
\int_{\Omega} \nabla\phi^\delta\cdot\mvec{u}^\delta \diff{V},
\end{eqnarray*}
for all test functions $\phi^\delta$ and $\mvec{w}^\delta$ in the
specified spaces.


\subsection{Properties of discretised Laplacian}

The first property to note for the \pdgp element pair is that the
discretised gradient in the $P1_{DG}$ velocity space $\mvec{q}^\delta$
of a P2 function $h^\delta$ obtained from
\[
\int_\Omega \mvec{w}^\delta\cdot\mvec{q}^\delta \diff{V}
= \int_\Omega \mvec{w}^\delta\cdot \nabla h^\delta \diff{V}
\]
for all $P1_{DG}$ test functions $\mvec{w}^\delta$, satisfies
\[
\q^\delta = \nabla h^\delta 
\]
at each point. To prove this, note that all continuous P2 functions $h$
have $P1_{DG}$ gradients. This means that we may choose
$\mvec{w}^\delta=\mvec{q}^\delta-\nabla h^\delta$ and hence
\[
\int_\Omega |\mvec{q}^\delta-\nabla h^\delta|^2
\diff{V}=0.
\]
Since $\mvec{q}^\delta$ and $\nabla h^\delta$ are piecewise
polynomials this means that they are identically equal.

The discretised Laplacian $L^\delta$ in the layer thickness space is
obtained by applying the discretised divergence to the discretised gradient
$\mvec{q}^\delta$:
\begin{eqnarray*}
\int_\Omega \phi^\delta L^\delta h^\delta \diff{V}&=&
-\int_\Omega \nabla\phi^\delta\cdot \mvec{q}^\delta \diff{V} 
+ \int_{\partial\Omega} \phi\vec{n}\cdot\mvec{q}^{\delta} \diff{S} \\
&=& 
-\int_\Omega \nabla\phi^\delta\cdot \nabla h^\delta \diff{V}
+ \int_{\partial\Omega} \phi{\pp{h}{n}}^{\delta} \diff{S},
\end{eqnarray*}
for any $P2$ test function $\phi^\delta$, which is the standard Galerkin
finite element discretisation of the Laplace operator obtained by
multiplying by a test function and integrating by parts.  The properties of
this operator are well-known in the finite element literature; in particular
it has no spurious eigenvalues.

\subsection{Exactly steady geostrophic modes}

In the linear shallow-water equations with rotation, the geostrophic
balanced modes with
\begin{equation}
\label{strong curl}
\mvec{u}_t=0 \implies \mvec{u} = \nabla^\perp \psi, 
\end{equation}
are steady in time, where
\[
\psi = \frac{\Ro}{\Fr^2}h, \quad \nabla^\perp \psi = (-\psi_y,\psi_x)
\]
 This is because $h_t=0$ since
$\nabla\cdot\mvec{u}=0$ for these modes. In this section we show that
the balanced states in discretisations with the \pdgp element pair are
also completely steady; this means that the \pdgp element pair
represents balanced states very well and so is ideal for shallow-water
ocean modelling.

The geostrophically balanced states in the finite element
discretisation satisfy
\begin{equation}
\label{discrete balance}
\dd{}{t}\int_\Omega \mvec{w}^\delta\cdot \mvec{u}^\delta\diff{V}
=0, \quad \implies
\int_\Omega \mvec{w}^\delta\cdot \mvec{u}^\delta\diff{V}
=\int_\Omega \mvec{w}^\delta \cdot \nabla^\perp \psi^\delta\diff{V},
\end{equation}
for all $P1^{DG}$ test functions $\mvec{w}^\delta$.

The property described in the previous section can be trivially
extended to show that the finite element velocity $\mvec{u}^\delta$
obtained from this equation for a given $\psi^\delta$ satisfies
equation \eqnref{strong curl} with $u=u^\delta$ and
$\psi=\psi^\delta$.  We can use this property to show that any
geostrophically balanced velocity field $\mvec{u}$ obtained from a
streamfunction $\psi$ which is constant on the boundary satisfies the
discrete divergence equation
\[
-\int_\Omega \nabla\phi^\delta\cdot\MM{u}^\delta\diff{V}=0,
\]
for all $P2$ test functions $\phi^\delta$.  To prove this, note that
\begin{eqnarray}
 -\int_\Omega \nabla \phi^\delta\cdot\mvec{u}^\delta\diff{V}
\label{weak curl1} & = & -\int_\Omega \nabla \phi^\delta\cdot\nabla^\perp 
\psi^\delta\diff{V} \\ 
& = & \sum_E\int_E \phi^\delta\underbrace{\nabla\cdot\nabla^\perp \psi^\delta
}_{=0}\diff{V} 
\label{bc1} 
- \sum_E\int_{\partial E}
\phi^\delta\mvec{n}\cdot\nabla^\perp \psi^\delta\diff{S} 
 \\
\label{jump}
& = & -\sum_\Gamma\int_\Gamma\phi^\delta 
\underbrace{[[\mvec{\mvec{n}}\cdot\nabla^\perp 
\psi^\delta]]}_{=0}\diff{V}- 
\int_{\partial\Omega}
\phi^\delta\underbrace{\mvec{u}\cdot\mvec{n}}_{=0}\diff{S} =0,
\end{eqnarray}
where $\sum_E$ indicates a sum over all elements $E$, $\partial E$ is the
boundary of element $E$, $\sum_\Gamma$ indicates a sum over all orientated
internal element boundaries in the mesh, and $[[f]]$ indicates the jump in a
function $f$ across a surface $\Gamma$. In equation \eqnref{bc1} the normal
component of velocity vanishes exactly on $\partial\Omega$ as $h$ is
constant on $\partial\Omega$ and the balanced velocity is obtained from the
pointwise curl of the streamfunction $\psi$. In \eqnref{jump} the jump in
the normal component of $\nabla^\perp \psi$ vanishes because the tangential
derivative of functions in $P2$ is continuous across element
boundaries.\footnote{The right-hand side of \eqnref{weak curl1} can be shown
  to vanish for general functions from the space $H^1$ (which contains the
  P2 functions) by taking a convergent sequence of smooth functions and
  passing to the limit. However, the extra property of continuous tangential
  derivatives of P$n$ functions facilitates the simpler proof given here.}

    The proof of this property is easily extended to the general
\pndgpn element pair 
\emph{i.e.}, $n$th order discontinuous velocity and 
$(n+1)$-th order continuous layer thickness. It is also easily extended to
the three-dimensional case in which $\mvec{u}=\nabla\wedge\mvec{\Psi}$
for any vector field $\mvec{\Psi}$ which is constant on the boundary.
\section{Numerical tests}
\label{numerical tests}
In this section we illustrate and explore the properties of the \pdgp
element applied to the linear rotating shallow-water equations.

\subsection{Representation of geostrophic balance}
\begin{figure}
\begin{center}
\includegraphics*[width=6cm]{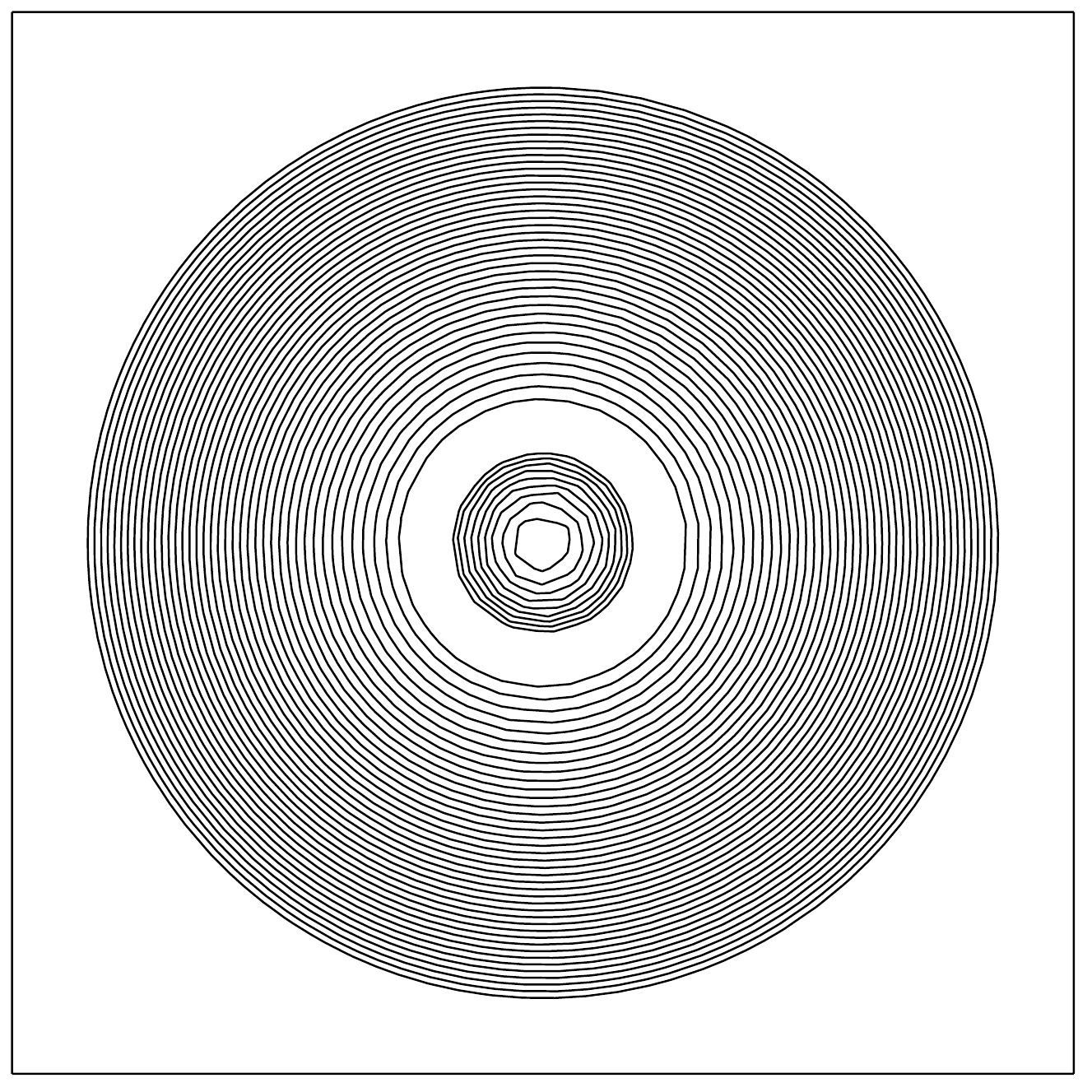}
\includegraphics*[width=6cm]{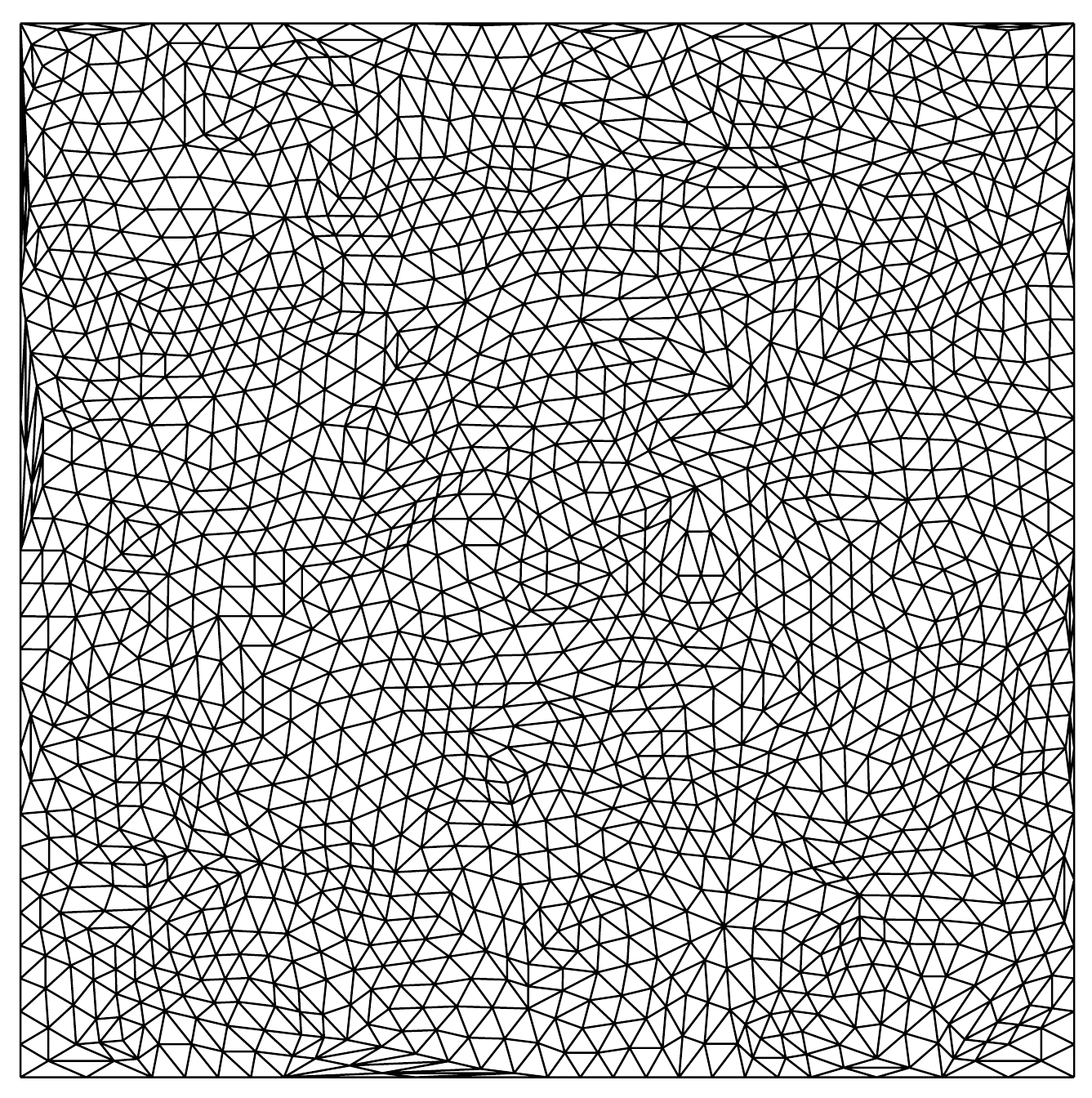}
\end{center}
\caption{\label{contours} {\bfseries Left}: Streamlines of balanced velocity
  obtained from a Gaussian streamfunction distribution using the \pdgp
  discretisation. The streamlines are very smooth showing that the
  discretisation does not introduce spurious oscillations. This is because
  in this case the balanced velocity can be obtained by taking the pointwise
  (strong) gradients of the streamfunction. {\bfseries Right}: The mesh used
  for this calculation. The mesh was deliberately distorted to illustrate
  that this property is not dependent on mesh quality.}
\end{figure}

\citet{Ro_etal1998} tested a number of element pairs for their
ability to represent geostrophic balance. This was done by selecting a
streamfunction field, computing the balanced velocity field from equation
\eqnref{discrete balance}, and plotting streamlines. Element pairs
were compared by the smoothness of the resulting streamlines on
structured and unstructured meshes.  Here we just note that, as
described in the previous section, the balanced velocity for the \pdgp
element is obtained from the pointwise gradients of streamfunction and so
streamlines of the discretised balanced velocity field are simply
contours of the discretised streamfunction field. This means that the
balanced velocity field is actually as accurate as possible for the
P2 streamfunction field. Plots of some resulting streamlines are given in
figure \ref{contours}; for comparison with other element pairs see
\citet{Ro_etal1998}.

\subsection{Steady states}
In \citet{Ro_etal1998}, another numerical test was performed in which
the linear rotating shallow-water equations were initialised in a
geostrophic state; streamlines were plotted after some time which
showed that the $\PP1$ $\iso$ $\PP2-\PP{0-3}$ element pair (proposed
in that paper) preserved the steady state to excellent accuracy. In
the case of the \pdgp element, we have already shown in the previous
section that geostrophic states are exactly steady so it remains to
verify this numerically. Using the mesh shown in figure
\ref{contours}, we computed randomly generated streamfunction fields
with $\psi=0$ on the boundary together with their geostrophically
balanced velocity fields obtained from equation \eqnref{strong curl},
and integrated the equations in time using the Crank-Nicholson
method. We observed that the layer thickness $h$ and velocity
$\mvec{u}$ remained constant up to machine precision, confirming that
the geostrophic modes are completely uncoupled from the
inertia-gravity waves. For the time evolution of geostrophic states
using other element pairs, see \citet{Ro_etal1998}.

\subsection{Kelvin waves}
\begin{figure}
\begin{center}
\includegraphics*[width=8cm]{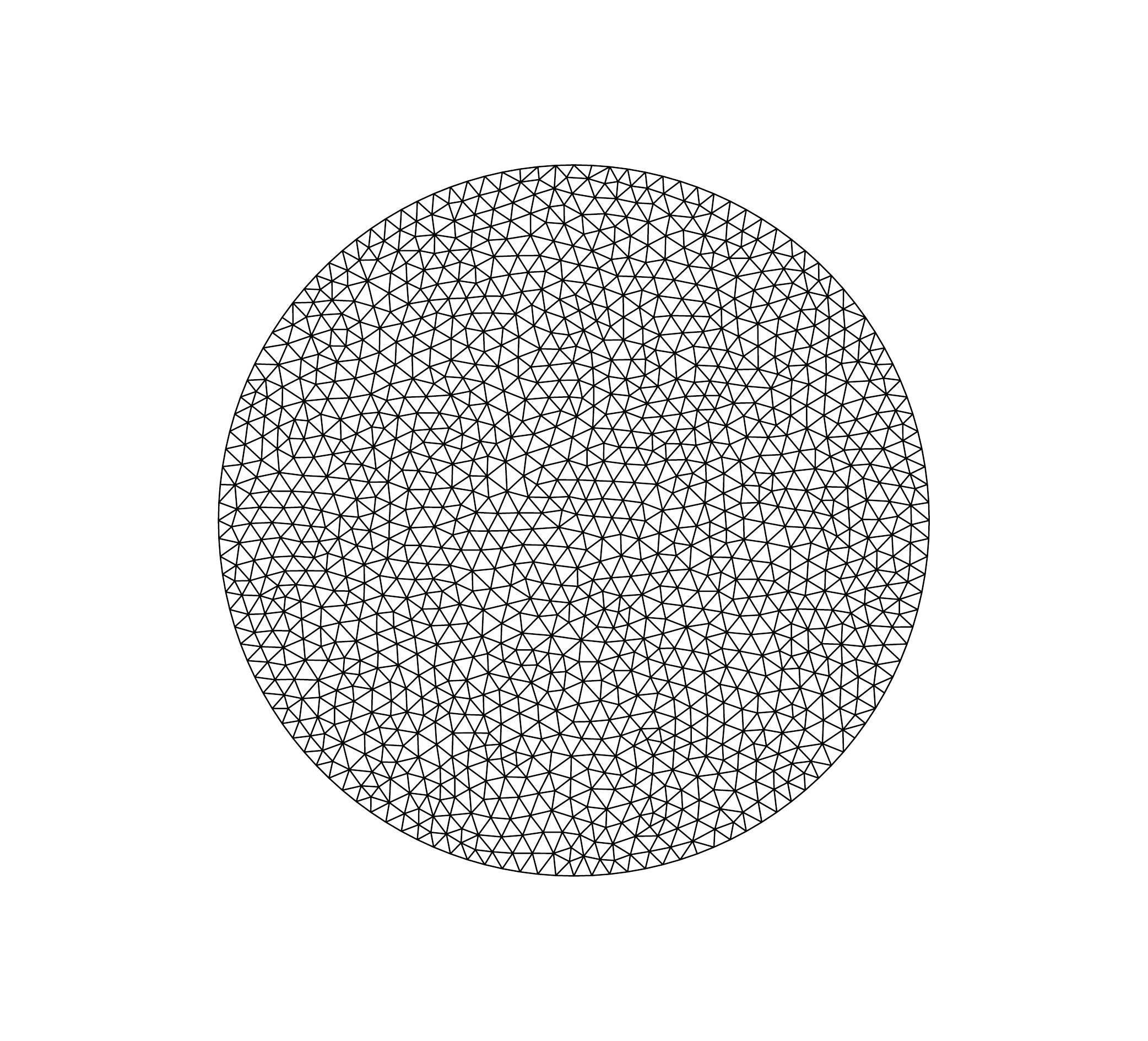}
\end{center}
\caption{\label{basinmesh}Plot showing the mesh used for Kelvin wave tests.}
\end{figure}

\begin{figure}
\begin{center}
\includegraphics*[width=4cm]{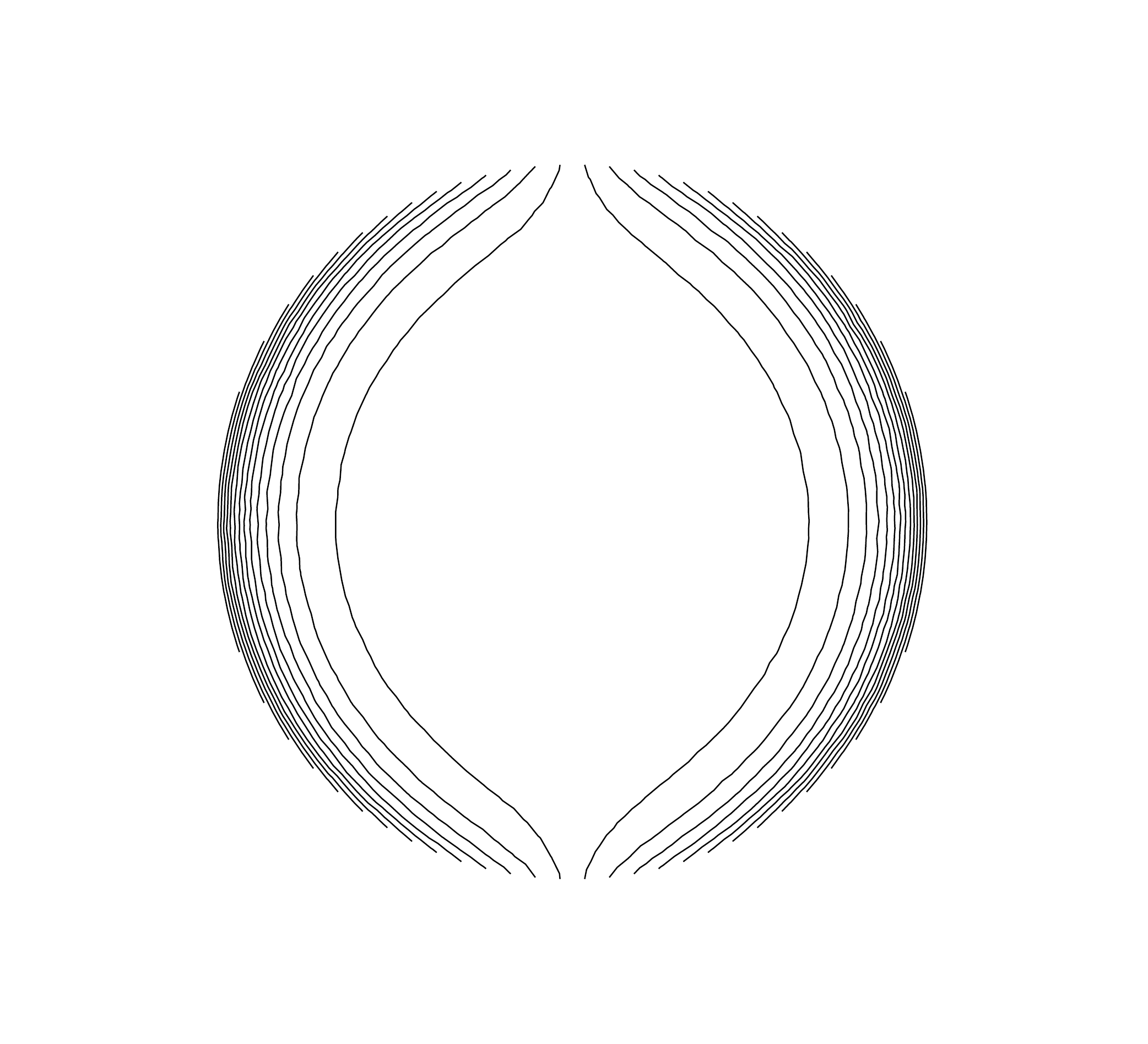}
\includegraphics*[width=4cm]{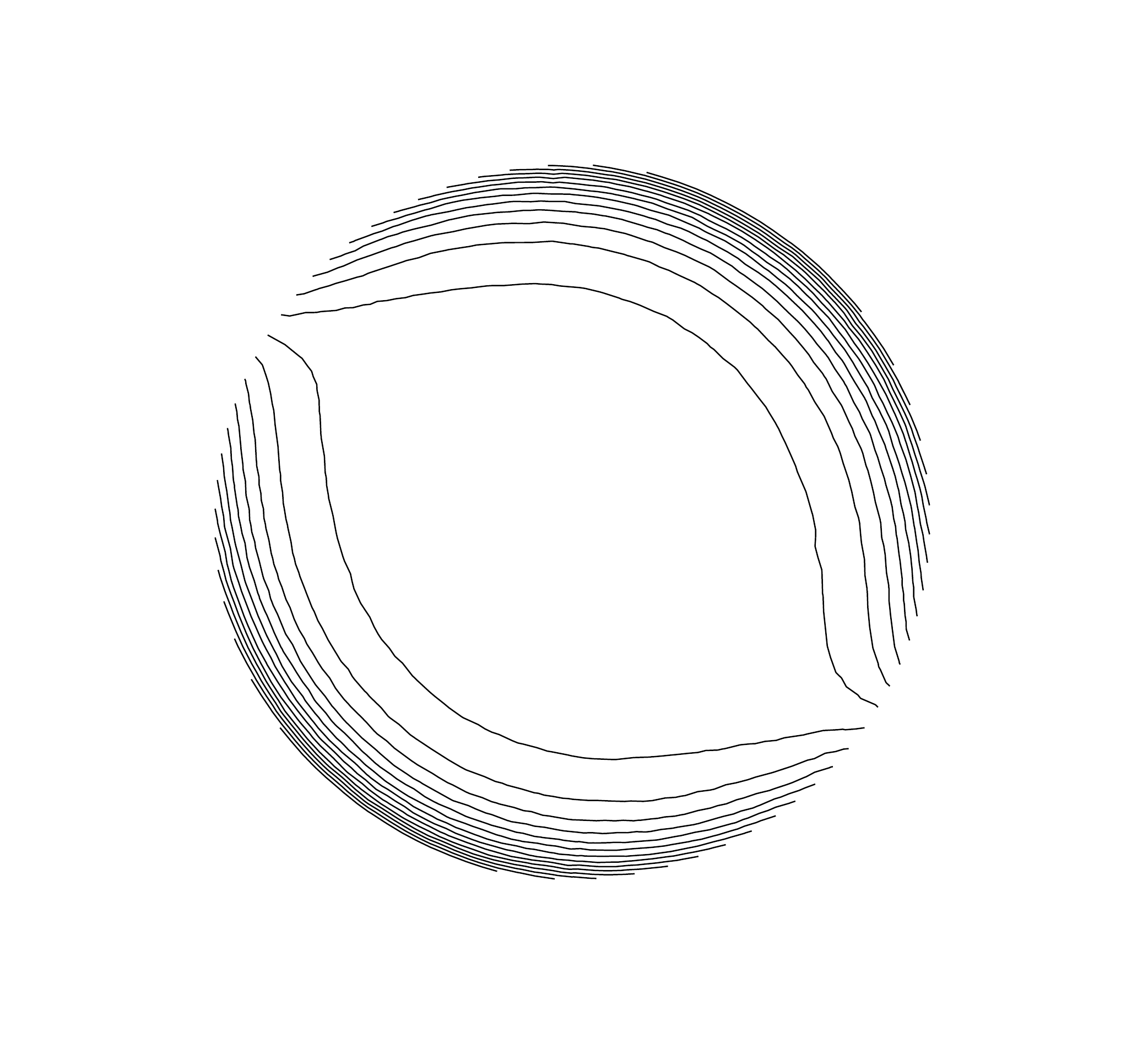} \\
\includegraphics*[width=4cm]{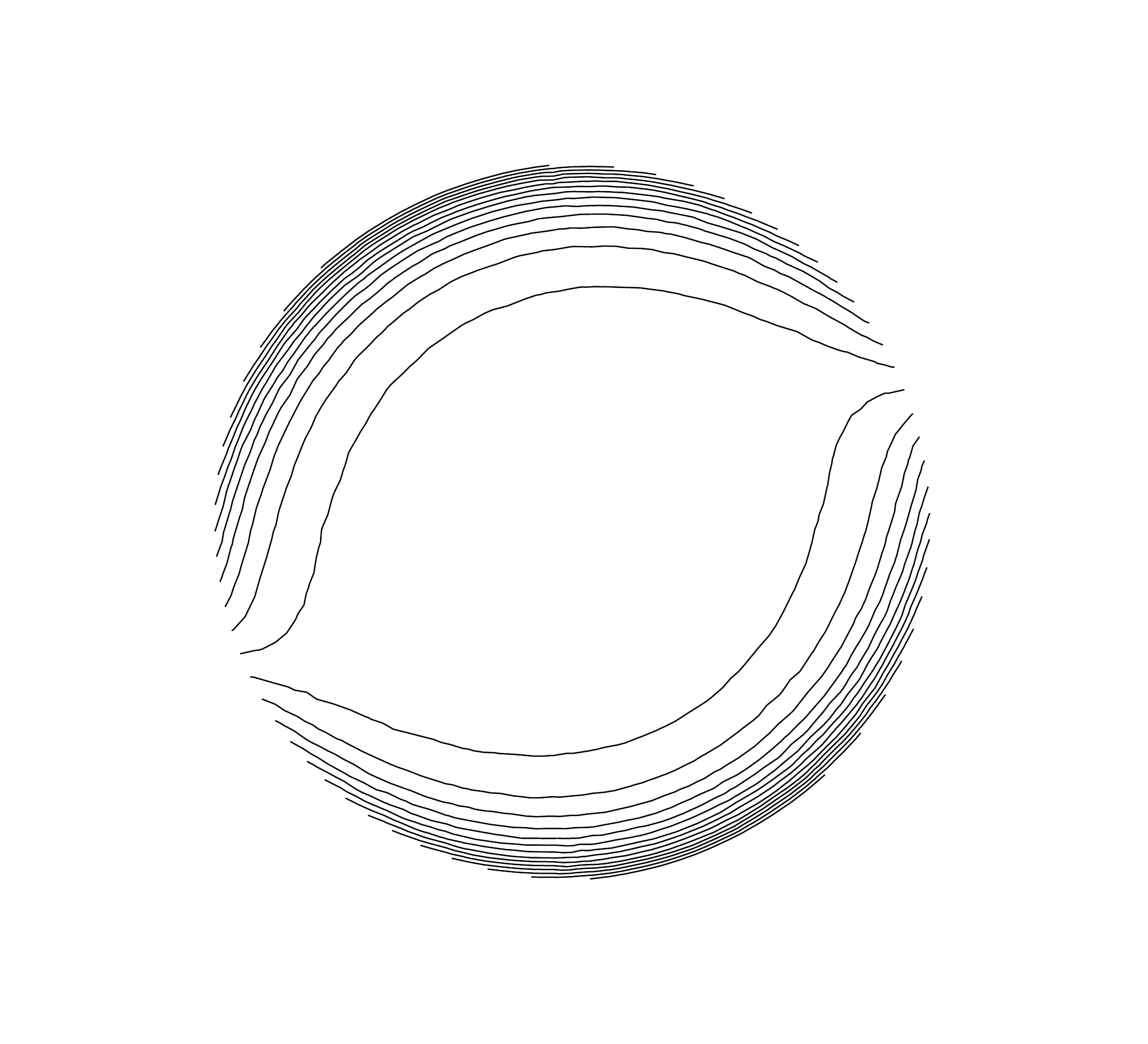}
\includegraphics*[width=4cm]{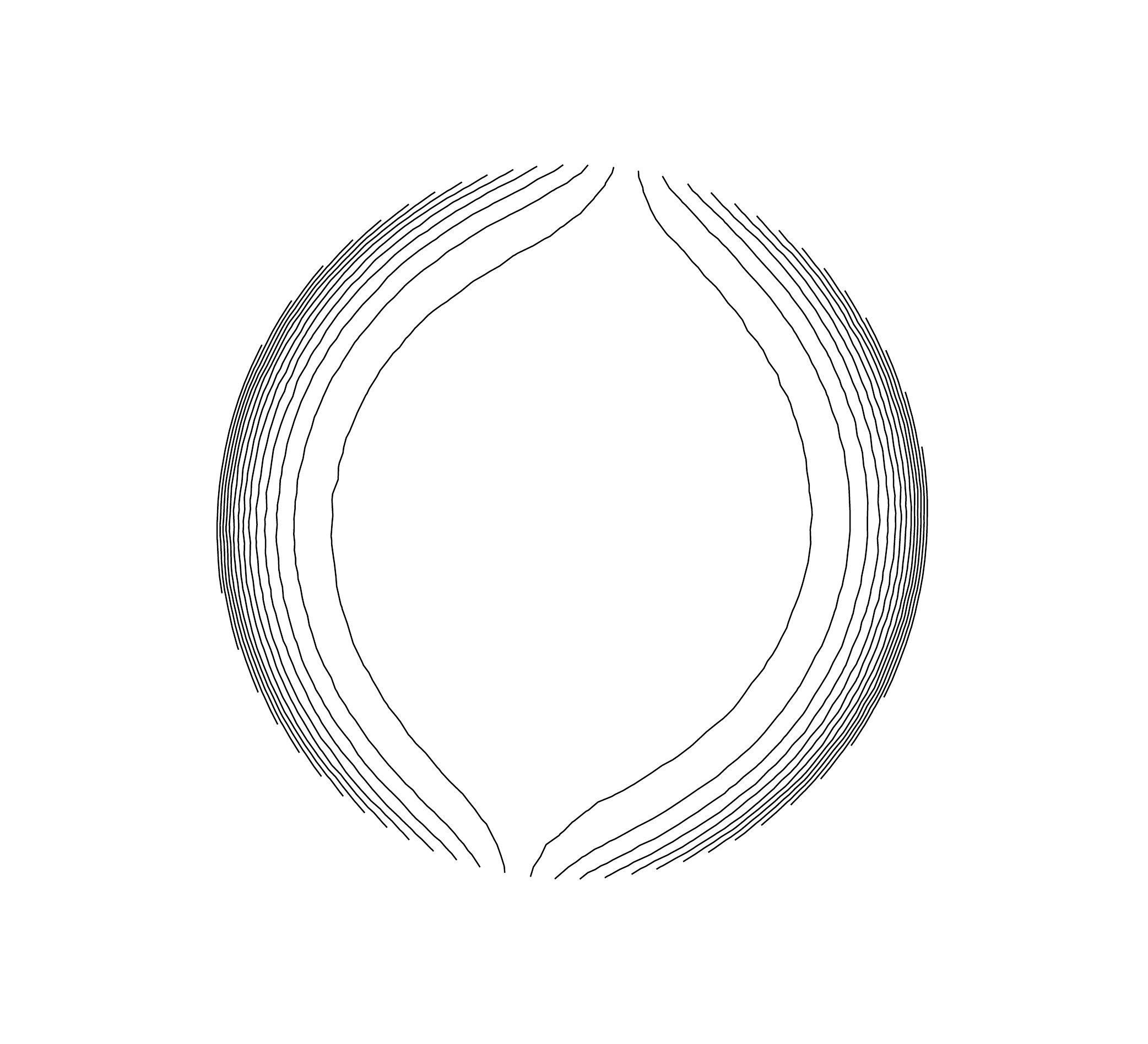}
\end{center}
\caption{\label{kelvin}Plots showing contours of layer thickness $h$
  at times $t=$0 (top left), 30 (top right), 60 (bottom left) and 90
  (bottom right) for the circular Kelvin wave test case. No spurious
  oscillations are observed, which verifies that the \pdgp element
  maintains geostrophic balance in the normal direction.}
\end{figure}

We tested the \pdgp element using a Kelvin wave initial condition; the
Kelvin wave is a trapped coastal wave which is geostrophically
balanced in the direction normal to the coast which propagates at the
fast gravity wave speed $1/\Fr$ for the case of a straight
coastline. The aim of the test is to verify that the Kelvin wave does
not shed any spurious inertia-gravity waves. We used the circular Kelvin
wave initial condition given by
\begin{eqnarray*}
h(r,\theta) & = & e^{(r-r_0)/\Ro}\cos\theta, \\
u_\theta(r,\theta) & = & \frac{1}{\Fr}e^{(r-r_0)/\Ro}\cos\theta, \\
u_r & = & 0,
\end{eqnarray*}
with $\Ro=0.1$ and $\Fr=1$. The Kelvin wave propagates around the circular
coast, maintaining geostrophic balance in the normal direction. The mesh
used for the discretisation is shown in figure \ref{basinmesh}. We
integrated the equations in time for $0>t>100$ using the Crank-Nicholson
method and a time step size $\Delta t=0.01$. Figure \ref{kelvin} shows the
layer thickness at various times: there are no spurious gravity waves
observed, which means that the \pdgp element pair is maintaining geostrophic
balance in the normal direction as well as the Kelvin wave structure.

\begin{figure}
\begin{center}
\includegraphics*[width=14cm]{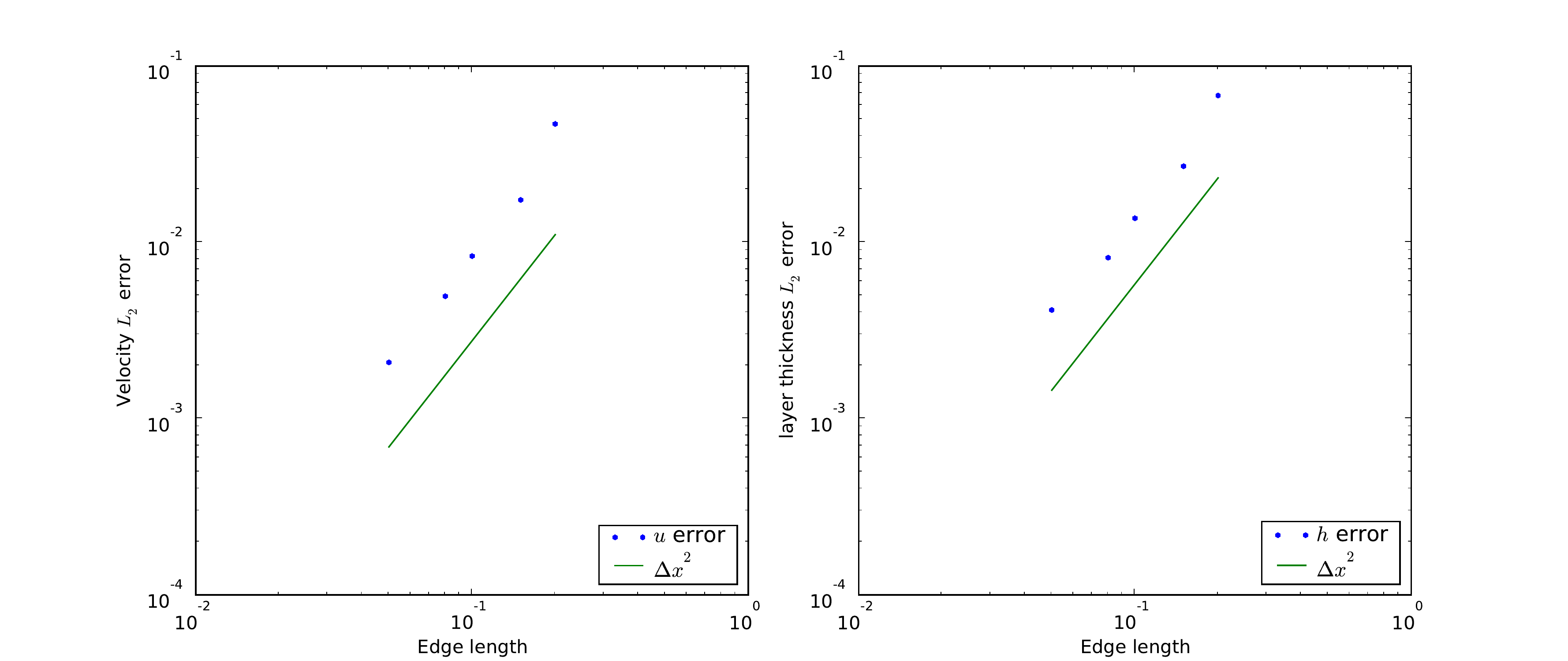}
\end{center}
\caption{\label{convergence} Convergence plots for tests with a Kelvin
  wave propagating along a flat coast performed on unstructured
  isotropic triangular meshes, showing error $\epsilon$ against
  element edge length $\Delta x$. {\bfseries Left}: $L_2$ error in
  velocity plotted against element edge length.  {\bfseries
    Right}: $L_2$ error in velocity plotted against element edge
  length.  Both plots show that the errors scale with $\Delta x^2$ as
  $\Delta x\to 0$.}
\end{figure}

To check convergence of the method we integrated a Kelvin wave in the
rectangular domain $\Omega=\{\mvec{x}:-15<x<15, 0<y<3\}$ with initial
condition
\[
h = e^{-y/\Ro}e^{-(x-5)^2}, \qquad 
\mvec{u} = (e^{-y/\Ro}e^{-(x-5)^2},0).
\]
If this initial condition is used in the domain $\Omega^\infty =
\{\mvec{x}:-\infty<x<\infty, 0<y<\infty\}$, then the equation has the
exact solution
\[
h = e^{-y/\Ro}e^{-(x+t/\Fr^2-5)^2}, \qquad 
\mvec{u} = (e^{-y/\Ro}e^{-(x+t/\Fr^2-5)^2},0).
\]
We integrated the system to time $t=10$. For this time interval the
solution is almost zero for $y>1$ and $|x|>6$ and so the exact
solution is a good approximation. The timestep was chosen to have a
wave Courant number of less than $0.1$ for all simulations so that the
errors are dominated by the spatial discretisation. We refined the
mesh isotropically in space in the region where the solution was
non-zero during the calculation and computed the $L_2$ errors of the
velocity and the layer thickness for various element edge lengths in
the refined region. Plots of the numerical errors are given in figure
\ref{convergence}.  A linear regression on these values showed that
the velocity errors were proportional to 2.19 and the layer thickness
errors were proportional to 1.98. These results suggest that the
errors in velocity and layer thickness in the spatial discretisation
scale quadratically with the edge length, as would be expected from
approximation theory. They are also an indication that the element
pair is stable: if the element pair were unstable then there would be
spurious modes present which would lead to slower convergence than
that expected from approximation theory.

\section{Summary and Outlook}
\label{summary}
In this paper we introduced the \pdgp element pair applied to the
linear shallow-water equations on an $f$-plane. We showed that the
element pair has the property that all geostrophically balanced states
which strongly satisfy the boundary conditions are exactly steady
since their discrete divergence is identically zero. This means that
the element pair has excellent geostrophic balance properties. We
verified these properties by computing the evolution of balanced
states, and by simulating Kelvin wave solutions which are
geostrophically balanced in one direction. Finally we gave convergence
test results which show that the numerical solutions have errors which
decay quadratically with element edge length; this verifies the
LBB-stability properties discussed in \citet{CoHaPaRe2007}.

In future work we shall compare this element pair with other low-order
element pairs such as the $P0_{DG}-P1$, $P1_{NC}-P1$ and $RT0$
pairs. Whilst the discontinuous velocity means that the \pdgp pair has
a large number of degrees of freedom per element, the remarkable
accuracy of the first half of the dispersion relation \citep[noted
in][]{CoHaPaRe2007} suggests that the element may be competitive,
especially given its excellent treatment of geostrophic balance, and
local conservation of momentum. The higher-order extensions such as
$\PP2_{DG}-\PP3$ will also be examined. We shall investigate the
performance of the element once nonlinear advection has been
introduced.

A key advantage of this element pair is that the extension to three
dimensions is also LBB-stable; the property that geostrophically
balanced states are exactly divergence-free also extends to the three
dimensional case. We shall investigate the performance of this element
pair in fully three-dimensional unstructured mesh ocean modelling in
the ICOM model \citep{PaPiGo2005}. We also expect that if the
buoyancy is discretised using $\PP1_{DG}$ elements, then the
discretisation will also preserve hydrostatic balance very well; this
will be investigated in future work.

\section{Acknowledgements}
We thank Greg Pavliotis for discussions about the proof that the
geostrophic modes are completely decoupled for this element pair, and
all of the AMCG team for their collaborative contributions. The
authors acknowledge funding from NERC consortium grant NE/C52101X/1.
\bibliography{DgCgCoriolis}

\begin{thebibliography}{12}
\expandafter\ifx\csname natexlab\endcsname\relax\def\natexlab#1{#1}\fi
\expandafter\ifx\csname url\endcsname\relax
  \def\url#1{\texttt{#1}}\fi
\expandafter\ifx\csname urlprefix\endcsname\relax\def\urlprefix{URL }\fi

\bibitem[{Ainsworth et~al.(2006)Ainsworth, Monk, and Muniz}]{AiMoMu2006}
Ainsworth, M., Monk, P., Muniz, W., 2006. Dispersive and dissipative properties
  of discontinuous {G}alerkin finite element methods for the second-order wave
  equation. J. Sci. Comput. 27~(1-3), 5--40.

\bibitem[{Ambati and Bokhove(2007)}]{AmBo2007}
Ambati, V., Bokhove, O., 2007. Space-time discontinuous galerkin discretization
  of rotating shallow water equations. Journal of Computational Physics
  225~(2), 1233--1261.

\bibitem[{Bernard et~al.(2007)Bernard, Chevaugeon, Legat, Deleersnijder, and
  Remacle}]{Be_etal2007}
Bernard, P.~E., Chevaugeon, N., Legat, V., Deleersnijder, E., Remacle, J.~F.,
  2007. High-order h-adaptive discontinuous galerkin methods for ocean
  modelling. Ocean Dynamics 57, 109--121.

\bibitem[{Cotter et~al.(2008)Cotter, Ham, Pain, and Reich}]{CoHaPaRe2007}
Cotter, C.~J., Ham, D.~A., Pain, C.~C., Reich, S., 2008. {LBB} stability of a
  mixed {G}alerkin finite element pair, submitted
  http://arxiv.org/0707.4607.

\bibitem[{Giraldo(2006)}]{Gi2006}
Giraldo, F.~X., 2006. High-order triangle-based discontinuous galerkin methods
  for hyperbolic equations on a rotating sphere. J. Comput. Phys. 214~(2),
  447--465.

\bibitem[{Gresho and Sani(2000)}]{GrSa2000}
Gresho, P.~M., Sani, R.~L., 2000. Incompressible Flow and the Finite Element
  Method, Volume 2, Isothermal Laminar Flow. Wiley.

\bibitem[{Karniadakis and Sherwin(2005)}]{KaSh2005}
Karniadakis, G. E.~M., Sherwin, S., 2005. Spectral/hp Element Methods for
  Computational Fluid Dynamics. Oxford Science Publications, Ch.~7.

\bibitem[{Le~Roux et~al.(1998)Le~Roux, Staniforth, and Lin}]{Ro_etal1998}
Le~Roux, D., Staniforth, A., Lin, C.~A., 1998. Finite elements for
  shallow-water equation ocean models. Monthly Weather Review 126~(7),
  1931--1951.

\bibitem[{Levin et~al.(2006)Levin, Iskandarani, and Haidvogel}]{Levin2006}
Levin, J., Iskandarani, M., Haidvogel, D., 2006. To continue or discontinue:
  Comparisons of continuous and discontinuous {G}alerkin formulations in a
  spectral element ocean model. Ocean Modelling 15, 56--70.

\bibitem[{Pain et~al.(2005)Pain, Piggott, Goddard, Fang, Gorman, Marshall,
  Eaton, Power, and de~Oliveira}]{PaPiGo2005}
Pain, C., Piggott, M., Goddard, A., Fang, F., Gorman, G., Marshall, D., Eaton,
  M., Power, P., de~Oliveira, C., 2005. Three-dimensional unstructured mesh
  ocean modelling. Ocean Modelling 10, 5--33.

\bibitem[{Raviart and Thomas(1977)}]{RaTh1977}
Raviart, Thomas, 1977. A mixed finite element method for 2nd order elliptic
  problems. In: Mathematical Aspects of the Finite Element Method. Lecture
  Notes in Mathematics. Springer, Berlin.

\bibitem[{Walters and Casulli(1998)}]{WaCa1998}
Walters, R., Casulli, V., 1998. A robust, finite element model for hydrostatic
  surface water flows. Communications in Numerical Methods in Engineering 14,
  931--940.

\end{thebibliography}



\end{document}